\documentclass[reqno,12pt, oneside]{amsart}

\usepackage{enumerate}
\usepackage[nobysame]{amsrefs}
\usepackage{enumitem}
\usepackage{amssymb}
\usepackage{amsmath}
\usepackage{amsthm}
\usepackage{mathtools}
\usepackage{amscd}
\usepackage{microtype}
\usepackage[right=2cm,left=2cm, top=2.5cm, bottom=2.5cm,includehead]{geometry}
\usepackage{amsfonts}
\usepackage{booktabs}
\usepackage{hyperref}

\usepackage{tgpagella}
\usepackage{eulervm}

\usepackage{chngcntr}
\counterwithin{table}{section}

\usepackage[backgroundcolor=lightgray]{todonotes}

\theoremstyle{plain}
\newtheorem{theorem}{Theorem}

\newtheorem{proposition}[theorem]{Proposition}

\theoremstyle{remark}
\newtheorem{remark}[theorem]{Remark}
\newtheorem{example}[theorem]{Example}

\theoremstyle{definition}
\newtheorem{definition}[theorem]{Definition}

\DeclarePairedDelimiter{\abs}{\lvert}{\rvert}

\DeclarePairedDelimiter{\norm}{\lVert}{\rVert}

\DeclarePairedDelimiterX{\dset}[2]{\lbrace}{\rbrace}{#1\;\delimsize|\;#2}

\DeclareMathOperator{\var}{var}

\newcommand{\ball}{{\overline{B}}}

\title{Josephy's theorem, revisited}

\author{Daria Bugajewska$^*$}
\address[D.~Bugajewska]{Faculty of Mathematics and Computer Science\\
 Adam Mickiewicz University in Pozna\'n\\
ul.\ Uniwersytetu Pozna\'nskiego 4\\
61-614 Pozna\'n\\
Poland}
\email[D.~Bugajewska]{dbw@amu.edu.pl}

\author{Piotr Kasprzak}
\address[P. Kasprzak]{Department of Nonlinear Analysis and Applied Topology\\
Faculty of Mathematics and Computer Science\\
Adam Mickiewicz University in Pozna\'n\\
ul.\ Uniwersytetu Pozna\'nskiego 4\\
61-614 Pozna\'n\\
Poland}
\email[P.~Kasprzak]{kasp@amu.edu.pl}


\date{\today\\
\hspace*{4mm}$^*$\emph{Corresponding author}: Daria Bugajewska (dbw@amu.edu.pl)}

\begin{document}
\begin{abstract}
The main goal of this note is to characterize the necessary and sufficient conditions for a composition operator to act between spaces of mappings of bounded Wiener variation in a normed-valued setting. The necessary and sufficient conditions for local boundedness of such operators are also discussed. 
\end{abstract}

\maketitle

\section{Introduction}

Composition operators (also known as autonomous superposition or Nemytskii operators) are one of the most basic, yet, at the same time, most important mappings of nonlinear analysis. The theory of such operators in various function or sequence spaces is very rich; we refer interested readers to the celebrated monograph by Appell and Zabrejko~\cite{AZ} for more details.

When we focus our attention at a specific function space, namely the space of real-valued functions of bounded Jordan variation, we will see that the central role is played here by the Josephy theorem. Proven in 1981 (see~\cite{J}) it says that the composition operator $C_f$, generated by $f\colon \mathbb R \to \mathbb R$, maps the space $BV_1([a,b],\mathbb R)$ into itself if and only if the function $f$ is locally Lipschitz continuous. In subsequent years, Josephy's result was extended to more general spaces of functions of bounded variation. For instance, Marcus and
Mitzel in~\cite{MM} dealt with the Riesz variation. Further, Ciemnoczolowski
and Orlicz in~\cite{CO} worked with the Young variation under the additional
assumption that the Young function $\phi$ and its inverse satisfy the
so-called $\Delta_2$-condition. Their result can be also used to
characterize composition operators acting in spaces of functions of
bounded Wiener variation by taking $\phi(t)= t^p$. Let us add that a similar fact for the generalized Young variation, where the function $\phi$ depends on a parameter, was
obtained in~\cite{BBL}. Finally, Pierce and Waterman extended Josephy's
theorem to the $\Lambda$-variation. (We refer readers interested in this
topic to the monograph on bounded variation by Appell et al.~\cite{ABM}; see
also~\cite{reinwand}.)

All the above-mentioned results, including those discussed in the
monographs~\cites{ABM, reinwand}, share one thing in common. Namely, they deal with composition operators acting in spaces of \emph{real-valued} functions of
bounded variation. Maps of bounded variation taking valves in an
arbitrary Banach space $E$ are considered in the book by Dudley and
Norvai\v{s}a~\cite{DN}. In particular, acting conditions for composition operators
between $BV_p([a,b],E)$ and $BV_q([a,b],E)$ are discussed in Section~6.5 of that work. Unfortunately, these conditions are only sufficient, and, as shown by Example~\ref{ex:not_locally_bdd} below, are not necessary.

The main goal of this note is to extend Josephy's result to an abstract
setting, that is, to prove a full characterization of the conditions under which
$C_f$ maps $BV_p([a,b],E)$ into $BV_q([a,b],E)$ (see Theorem 4); here $E$ is a normed space and $1\leq p \leq q<+\infty$. Notably, in contrast to the real-valued case, these necessary and sufficient acting conditions do not guarantee the local boundedness of the composition operator $C_f$ (see Example~\ref{ex:not_locally_bdd}). Therefore, in Theorem~\ref{thm:loc_bounded} we present necessary and sufficient conditions under which the considered composition operator is locally bounded.

\section{Notation}

Throughout the paper by $E$ we will always denote a normed space (over either the field of real or complex numbers) endowed with the norm $\norm{\cdot}_E$. Given a point $x \in E$ and a positive number $r$, the symbol $\ball_E(x,r)$ will stand for the closed ball in $E$ with center $x$ and radius $r$.

In the paper we will also use two notions related to compactness. A non-empty set $A \subseteq E$ is called \emph{relatively compact}, if its closure with respect to $E$ is compact. Equivalently, $A\subseteq E$ is relatively compact, if each sequence of elements of $A$ has a subsequence convergent to a point in $E$. On the other hand, a non-empty set $A\subseteq E$ is called \emph{precompact}, if each sequence of its elements contains a Cauchy subsequence. From the definitions it is obvious that precompactness is a notion weaker than relative compactness, and that in complete normed spaces those two notions coincide.

When talking about sums of some elements, we will also adopt the convention that if the upper summation index is greater than the lower one, then the sum under consideration, by definition, equals to zero. 

\section{Composition operators}

The main object of our study are composition operators acting on some function spaces.
So let us proceed with their definition. Let $E$ be a normed space and by $X$ let us denote a set of some $E$-valued maps defined on a compact interval $[a,b]$.
A \emph{composition operator} (or, \emph{autonomous superposition operator}) $C_f$,
generated by $f \colon E \to E$, is a map which to each $x \in X$ assigns the function $C_f(x)$ given by $C_f(x)(t)=f(x(t))$ for $t \in [a,b]$. In general, it
may happen that $C_f(x)$ does not belong to $X$, and, of course, the most interesting situation is when it does. 

In this paper $X$ will be the space of maps of bounded Wiener variation, and it will turn out that the most suitable choice of generators $f$ will revolve around classes of H\"older continuous maps. Let us recall that $f\colon E\to E$ is
said to be \emph{H\"older continuous} on precompact (resp., bounded) subsets of $E$ with exponent $\alpha \in (0,1]$, if for every precompact (resp., bounded) set $K \subseteq E$ there exists a constant $L_K\geq 0$ such that $\norm{f(u)-f(w)}_E \leq L_K \norm{u-w}_E^\alpha$ for any $u,w \in K$.

\section{Maps of bounded Wiener variation}

The aim of this section is to recall the definition of mappings of bounded Wiener variation and present some of their properties. 

\begin{definition}
Assume that $E$ is a normed space and that $p \in [1,+\infty)$. Let $[a,b]$ be a non-empty interval, $[c,d] \subseteq [a,b]$ and let $x \colon [a,b] \to E$ be a map. Then, the (possibly infinite) quantity
\[
 \var_p(x, [c,d]) = \sup \sum_{i=1}^n \norm[\big]{x(t_i)-x(t_{i-1})}_E^p,
\] 
where the supremum is taken over all finite partitions $c=t_0<t_1<\ldots<t_n=d$ of the interval $[c,d]$, is called the \emph{$p$-variation} \textup(or, \emph{the Wiener variation}\textup) of the map $x$ over the interval $[c,d]$. If $\var_p(x,[a,b])<+\infty$, we say that $x$ is of \emph{bounded $p$-variation} (or just a \emph{$BV_p$-function}).
\end{definition}

Very often, when the interval $[a,b]$, on which the mapping $x$ is defined, will be known from the context, instead of $\var_p(x,[a,b])$ we will simply write $\var_p x$. As in the real-valued setting, it is easy to check that the set $BV_p([a,b],E)$ of all maps $x \colon [a,b] \to E$ of bounded $p$-variation, when endowed with the norm $\norm{x}_p:=\norm{x(a)}_E+(\var_p x)^{1/p}$, is a normed space. It can be also shown that $(\var_q x)^{1/q} \leq (\var_p x)^{1/p}$ for $x \in BV_p([a,b],E)$, when $1\leq p \leq q < +\infty$ (cf.~\cite{ABM}*{Proposition~1.38}). Hence, $BV_p([a,b],E)\subseteq BV_q([a,b],E)$ if $1\leq p \leq q < +\infty$.

From the definition of a mapping of bounded Wiener variation it easily follows that each such a map must be bounded. Actually, it can be shown that $\norm{x}_{\infty} \leq \norm{x}_p$ for every $x \in BV_p([a,b],E)$; here, by $\norm{\cdot}_\infty$ we denote the supremum norm of $x$, that is, $\norm{x}_\infty:=\sup_{t \in[a,b]}\norm{x(t)}_E$. However, an even stronger result is true as evidenced by the following proposition. (Since we could not find any book containing this fact, we decided to provide its proof.)

\begin{proposition}\label{prop:precompact_range}
Maps of bounded $p$-variation for $p \geq 1$ have precompact ranges.
\end{proposition}

\begin{proof}
Let $x \colon [a,b] \to E$ be a map of bounded $p$-variation. Consider a sequence $(x(t_n))_{n \in \mathbb N}$ of elements of the range of $x$. Passing to a subsequence if necessary, we may assume that the sequence $(t_n)_{n \in \mathbb N}$ is convergent to a real number $t \in [a,b]$. If there are infinitely many indices $n \in \mathbb N$ such that $t=t_n$, then the sequence $(x(t_n))_{n \in \mathbb N}$ admits a constant (and hence a Cauchy) subsequence, and the proof is complete. Therefore, without loss of generality, we may assume that $t\neq t_n$ for all $n \in \mathbb N$. Now, we can extract from $(t_n)_{n \in \mathbb N}$ a strictly monotone subsequence. By abuse of notation we will denote it using the same symbol $(t_n)_{n \in \mathbb N}$, and moreover, we will assume that it is increasing. (The case when $(t_n)_{n \in \mathbb N}$ is decreasing can be treated in a similar manner.) We claim that $(x(t_n))_{n \in \mathbb N}$ is Cauchy. Suppose that this is not the case. Then, there exists a positive constant $\varepsilon>0$ and two subsequences $(\tau_n)_{n \in \mathbb N}$, $(\sigma_n)_{n \in \mathbb N}$ of $(t_n)_{n \in \mathbb N}$ such that $a<\tau_1<\sigma_1<\tau_2<\sigma_2<\ldots<t\leq b$ and $\norm{x(\tau_n)-x(\sigma_n)}_E\geq \varepsilon$ for $n \in \mathbb N$. And so,
\[ 
 \var_p x \geq \sup_{m \in \mathbb N} \sum_{n=1}^m \norm[\big]{x(\tau_n)-x(\sigma_n)}_E^p =+\infty.
\]
This contradicts the fact that $x \in BV_p([a,b],E)$. Consequently, the sequence $(x(t_n))_{n \in \mathbb N}$ is Cauchy as claimed, and the range $x([a,b])$ is precompact.
\end{proof}

The following example shows that, in general, the range of a $BV_p$-map is not relatively compact.

\begin{example}\label{ex:not_rc}
Consider the space $c_{00}(\mathbb R)$ consisting of all real sequences with only finitely many non-zero terms. It is well-known that endowed with the supremum norm $\norm{(\xi_n)_{n \in \mathbb N}}_\infty:=\sup_{n \in \mathbb N}\abs{\xi_n}$, $c_{00}(\mathbb R)$ is not complete. In particular, the classes of its relatively compact and precompact subsets do not coincide. Let us consider the map $x \colon [0,1] \to c_{00}(\mathbb R)$ defined by the formula
\[
 x(t)=\begin{cases}
       (1,\frac{1}{2^2},\ldots,\frac{1}{k^2},0,0,0,\ldots), & \text{if $t \in [1-\frac{1}{k},1-\frac{1}{k+1})$ for some $k \in \mathbb N$,}\\
			 (0,0,0,\ldots), & \text{if $t=1$.}
			\end{cases}
\] 
Take an arbitrary finite partition $0=t_0<t_1<\ldots<t_n=1$ of the interval $[0,1]$. As adding new points to the partition may only increase the sum $\sum_{i=1}^n \norm{x(t_i)-x(t_{i-1})}_\infty$, we may assume that $t_{n-1}=1-\frac{1}{k+1}$ for some $k \geq 2$ and that all the points $\frac{1}{2},\frac{1}{3},\ldots,1-\frac{1}{k}$ appear among the points $t_1,\ldots,t_{n-2}$. Then,
\begin{align*}
 \sum_{j=1}^n \norm[\big]{x(t_j)-x(t_{j-1})}_\infty &\leq \sum_{i=1}^k \norm[\big]{x(1-\tfrac{1}{i+1})-x(1-\tfrac{1}{i})}_\infty + \norm[\big]{x(1-\tfrac{1}{k+1})}_\infty\\
& \leq 1 + \sum_{i=1}^\infty \frac{1}{(i+1)^2}<+\infty.
\end{align*}
Therefore, $x \in BV_1([0,1],c_{00}(\mathbb R))$. However, it is clear that the range $x([0,1])=\{(0,0,0,\ldots)\} \cup\dset{(1,\frac{1}{2^2},\ldots,\frac{1}{n^2},0,0,0,\ldots)}{n \in \mathbb N}$ is not a relatively compact subset of $c_{00}(\mathbb R)$.
\end{example}

\section{Josephy's theorem, revisited}

A composition operator $C_f$ acting on a function space $X$ assigns to 
each element $x \in X$ a new map $C_f(x)$. In general, this new mapping does not necessarily belong to $X$. Therefore, the most fundamental problem in the theory of composition operators is to provide necessary and sufficient conditions for $C_f$ to map the given space $X$ into itself.

In the case of the space $BV_1([a,b],\mathbb R)$ of real-valued functions of bounded
Jordan variation such conditions were given by Josephy in~\cite{J}. Here, we want to extend that result to maps with values in an arbitrary normed space. While we will follow Josephy's general approach, it is worth noting that our proof will not be a mere rewriting of the original argument. Firstly, we will have to be more subtle when dealing with certain sequences, since instead of (relative) compactness we will relay on a weaker notion -- precompactness. Secondly, showing that a function constructed in the proof is of bounded $p$-variation will require more effort. In contrast to the Jordan variation (or $p$-variation for $p=1$), adding new points to the given partition of $[a, b]$ could result in the variation sum decreasing. To see this phenomenon it suffices to consider the function $x \colon [0, 3] \to \mathbb R$ given by
\[
 x(t)=\begin{cases}
       0, & \text{if $t \in [0,1)$,}\\
			 1, & \text{if $t \in [1,2)$,}\\
			 2, & \text{if $t \in [2,3]$.}
			\end{cases}
\]
Then, $\abs{x(3)-x(0)}^2=4$, but $\abs{x(1)-x(0)}^2 + \abs{x(2)-x(1)}^2 + \abs{x(3)-x(2)}^2=2$. Last but not least, we will work with composition operators acting between spaces of bounded Wiener variation for different values of the parameter $p$.

After this lengthy introduction let us proceed to the main result of our paper.

\begin{theorem}\label{thm:jospehy}
Let $E$ be a normed space and let $1 \leq p \leq q < +\infty$. The composition operator $C_f$, generated by the mapping $f \colon E \to E$, maps the space $BV_p([a,b], E)$ into $BV_q([a,b],E)$ if and only if $f$ is H\"older continuous on precompact subsets of $E$ with exponent $p/q$. 
\end{theorem}

\begin{proof}
We begin with the proof of the sufficiency. Let $x \in BV_p([a,b], E)$. In view of Proposition~\ref{prop:precompact_range}, $x$ has precompact range. Therefore, setting $K:=x([a,b])$ there exists a constant $L_K\geq 0$ such that $\norm{f(x(t))-f(x(s))}_E \leq L_K \norm{x(t)-x(s)}^{p/q}_E$ for $t,s \in [a,b]$. From this it easily follows that $\var_q C_f(x) \leq L_K^q \var_p x$. In particular, $C_f(x) \in BV_q([a,b],E)$.

Now, we move to the necessity part. It will be divided into several steps.

\begin{enumerate}[label={\emph{Step \arabic*}.}, wide]
 \item \emph{Reduction}. It is fairly obvious that $f$ generates the composition operator acting between $BV_p([a,b], E)$ and $BV_q([a,b],E)$ if and only if it generates the composition operator acting between $BV_p([0,1], E)$ and $BV_q([0,1],E)$. Therefore, in our discussion we will assume that $[a,b]=[0,1]$.

\item \emph{The map $f$ is bounded on precompact subsets of $E$}. Suppose that there exist a precompact subset $A$ of $E$ and a sequence $(v_n)_{n \in \mathbb N}$ in $A$ such that $\norm{f(v_{n+1})}_E \geq \norm{f(v_n)}_E + 1$ for $n \in \mathbb N$. Passing to a subsequence, we may assume that $(v_n)_{n \in \mathbb N}$ is Cauchy and $\norm{v_{n+1}-v_n}_E \leq 1/n^2$ for $n \in \mathbb N$. Define $y \colon [0,1] \to E$ by the formula
\[
 y(t)=\begin{cases}
      v_n, & \text{if $t \in [1-\frac{1}{n},1-\frac{1}{n+1})$ for some $n \in \mathbb N$,}\\
			 0, & \text{if $t=1$.}
			\end{cases}
\] 
Using a similar reasoning to the one we used in Example~\ref{ex:not_rc} it can be shown that $y \in BV_1([0,1],E)\subseteq BV_p([0,1],E)$. However, for any $m \in \mathbb N$ we have
\begin{align*}
 \var_q C_f(y) \geq \sum_{n=1}^m \norm[\big]{f(y(1-\tfrac{1}{n+1}))-f(y(1-\tfrac{1}{n}))}_E^q = \sum_{n=1}^m \norm[\big]{f(v_{n+1})-f(v_n)}_E^q \geq m.
\end{align*}
Hence, $\var_q C_f(y)=+\infty$ -- a contradiction. This means that $f$ is bounded on every precompact subset of $E$.

\item \emph{Setting}. Suppose that the generator $f$ is not H\"older continuous with exponent $p/q$ on a precompact set $K\subseteq \ball_E(0,r)$. Then, for every $n \in \mathbb N$ we can find points $u_n,w_n \in K$ such that
\[
 \norm[\big]{f(u_n)-f(w_n)}_E > 4Mn^2 \norm[\big]{u_n-w_n}_E^{p/q},
\]
where $M:=\sup_{v \in K} \norm{f(v)}_E<+\infty$. Let us follow with some observations.
\begin{enumerate}[label=\textup{(\alph*)}]
 \item The points $u_n$ and $w_n$ must be distinct for every $n \in \mathbb N$. 

  \item For every $n \in \mathbb N$ we have
	\[
	 \norm[\big]{u_n-w_n}_E\leq \norm[\big]{u_n-w_n}_E^{p/q} < \frac{1}{4Mn^2} \cdot  \norm[\big]{f(u_n)-f(w_n)}_E \leq \frac{2M}{4Mn^2}=\frac{1}{2n^2}.
	\]
 \item Passing to subsequences if necessary, we may assume that $(u_n)_{n \in \mathbb N}$, $(w_n)_{n \in \mathbb N}$ are Cauchy and that $\norm{u_{n}-u_k}_E\leq 1/k^2$ and $\norm{w_{n}-w_k}_E\leq 1/k^2$ for $n \in \mathbb N$ and $1\leq k \leq n$.
\end{enumerate}

\item \emph{Definition of the map $x$}. For $n \in \mathbb N$ let $m_n:=[n^{-2q}\norm{u_n-w_n}_E^{-p}]$, where $[t]$ denotes the greatest integer smaller than or equal to $t$. Note that $m_n$ is well-defined and $m_n\geq 2$. Now, choose $m_n$ arbitrary points $1/(n+1)<a_1^n<a_2^n<\ldots<a_{m_n}^n<1/n$ and define $x \colon [0,1] \to E$ by the formula
\[
 x(t)=\begin{cases}
       0, & \text{if $t=0$,}\\
			 u_n, & \text{if $t=a_k^n$ for some $n \in \mathbb N$ and $k=1,\ldots,m_n$,}\\
			 w_n, & \text{if $t \in [\frac{1}{n+1},\frac{1}{n})\setminus\{a_1^n,\ldots,a_{m_n}^n\}$ for some $n \in \mathbb N$,}\\
			 w_1, & \text{if $t=1$.}
      \end{cases}
\]

\item \emph{The map $x$ is of bounded $p$-variation}. We claim that for every $n \in \mathbb N$ we have
\begin{equation}\label{eq:estimate_S}
 \var_p(x,[\tfrac{1}{n+1},1]) \leq \sum_{i=1}^{n} \frac{2}{i^2} + \sum_{i=1}^{n-1} \frac{2^{p+1}}{i^{2p}}.
\end{equation}
The estimate is clearly true for $n=1$, as in this case we have $\var_p(x,[\frac{1}{2},1]) \leq 2m_1\norm{u_1-w_1}_E^p \leq 2$. Assume now that~\eqref{eq:estimate_S} holds for all natural indices smaller than or equal to some $n \in \mathbb N$ and consider an arbitrary finite partition $1/(n+2)=\tau_0<\tau_1<\ldots<\tau_k=1$ of $[\frac{1}{n+2},1]$. Let $\tau_j$ be the greatest point of this partition that $1/(n+2)\leq \tau_j < 1/(n+1)$. Then, the next point $\tau_{j+1}$ lies in some interval $[\frac{1}{l+1},\frac{1}{l})$ for some $l \in \{1,\ldots,n\}$, or $\tau_{j+1}=\tau_k=1$. Hence, 
\[
 \sum_{i=1}^j \norm[\big]{x(\tau_i)-x(\tau_{i-1})}_E^p \leq 2m_{n+1}\norm[\big]{u_{n+1}-w_{n+1}}^p_E \leq \frac{2}{(n+1)^2}.
\]
Also, 
\[
 \sum_{i=j+2}^k \norm[\big]{x(\tau_i)-x(\tau_{i-1})}_E^p \leq \var_p(x,[\tfrac{1}{l+1},1]) \leq \sum_{i=1}^{l} \frac{2}{i^2} + \sum_{i=1}^{l-1} \frac{2^{p+1}}{i^{2p}}.
\]
Note that $x(\tau_j) \in \{u_{n+1},w_{n+1}\}$ and $x(\tau_{j+1})\in \{u_l,w_l\}$. So, $\norm{x(\tau_{j+1})-x(\tau_{j})}_E^p \leq 2^{p+1} l^{-2p}$. All of this implies that
\begin{align*}
 \sum_{i=1}^j \norm[\big]{x(\tau_i)-x(\tau_{i-1})}_E^p + \norm[\big]{x(\tau_{j+1})-x(\tau_{j})}_E^p + \sum_{i=j+2}^k \norm[\big]{x(\tau_i)-x(\tau_{i-1})}_E^p \leq \sum_{i=1}^{n+1} \frac{2}{i^2} + \sum_{i=1}^{n} \frac{2^{p+1}}{i^{2p}}.
\end{align*} 
Thus,
\[
  \var_p(x,[\tfrac{1}{n+2},1]) \leq \sum_{i=1}^{n+1} \frac{2}{i^2} + \sum_{i=1}^{n} \frac{2^{p+1}}{i^{2p}}
\]
and our claim follows by mathematical induction.

To estimate the $p$-variation of $x$ over the interval $[0,1]$ let  $0=t_0<t_1<\ldots<t_\lambda=1$ be an arbitrary finite partition of the interval $[0,1]$. Also, let $n \in \mathbb N$ be so that $1/(n+1)<t_1$. Then,
\begin{align*}
 &\sum_{i=1}^\lambda \norm[\big]{x(t_i)-x(t_{i-1})}_E^p\\
 &\qquad \leq 2^p\norm[\big]{x(\tfrac{1}{n+1})-x(0)}_E^p + 2^p\norm[\big]{x(\tfrac{1}{n+1})-x(t_1)}_E^p+ \sum_{i=2}^\lambda \norm[\big]{x(t_i)-x(t_{i-1})}_E^p\\
&\qquad \leq 2^p r^p + \sum_{i=1}^{\infty} \frac{2^{p+1}}{i^2} + \sum_{i=1}^{\infty} \frac{2^{2p+1}}{i^{2p}}<+\infty.
\end{align*}  
This shows that $\var_p(x,[0,1]) <+\infty$.

\item \emph{The map $C_f(x)$ is not of bounded $q$-variation}. In each interval $(\frac{1}{n+1},\frac{1}{n})$ let us choose $m_n$ points $b_1^n, \ldots,b_{m_n}^n$ so that $1/(n+1)<b_1^n<a_1^n<b_2^n<a_2^n<\ldots<b_{m_n}^n<a_{m_n}^n<1/n$. Then, for every $n \in \mathbb N$ we have
\begin{align*}
 \var_q C_f(x) &\geq \sum_{j=1}^n \sum_{i=1}^{m_j} \norm[\big]{f(x(b_i^j))-f(x(a_i^j))}^q_E =\sum_{j=1}^n {m_j} \norm[\big]{f(u_j)-f(w_j)}^q_E\\
 &\geq \sum_{j=1}^n 2^q m_j M^q j^{2q} \norm[\big]{u_j-w_j}^p_E \geq \sum_{j=1}^n M^q = nM^q. 
\end{align*}
Hence, $\var_q C_f(x)=+\infty$. This, however, contradicts the assumption that $C_f$ maps $BV_p([0,1],E)$ into $BV_q([0,1],E)$. Therefore, $f$ is H\"older continuous on precompact subsets of $E$ with exponent $p/q$. \qedhere
\end{enumerate}
\end{proof}

In the context of real-valued functions, one of the corollaries to Josephy-type theorems states that each composition operator acting between $BV_p([a,b],\mathbb R)$ and $BV_q([a,b],\mathbb R)$ for $1\leq p \leq q<+\infty$ maps bounded subsets of the former space into bounded subsets of the latter. (This property is commonly referred to as the \emph{local boundedness} of an operator; note that we do not assume here that $C_f$ is continuous.) However, in the abstract setting, such a result is no longer true, as demonstrated by the following example. 

\begin{example}\label{ex:not_locally_bdd}
Let $l^2(\mathbb R)$ be the Banach space of all real sequences that are square summable, endowed with the norm $\norm{(\xi_n)_{n \in \mathbb N}}_{l^2} = \bigl(\sum_{n=1}^\infty \abs{\xi_n}^2 \bigr)^{1/2}$. Let us also consider the map $f \colon l^2(\mathbb R) \to l^2(\mathbb R)$ given by the formula $f((\xi_n)_{n \in \mathbb N})=(\sup_{n \in \mathbb N} n(2\abs{\xi_n}-1), 0, 0, \ldots)$. It can be shown that $f$ is Lipschitz continuous on relatively compact (and hence precompact) subsets of $l^2(\mathbb R)$ -- cf.~\cite{CLM}*{pp.~63--64} and~\cite{cobzas_book}*{Theorem~2.1.6 and Remark~2.1.7}. Hence, in view of Theorem~\ref{thm:jospehy}, it generates the composition operator $C_f$ that maps $BV_1([0,1],l^2(\mathbb R))$ into itself. However, this operator is not locally bounded. To see this by $e_k$ let us denote the $k$-th unit vector of $l^2(\mathbb R)$, that is, the sequence consisting of all zeros and a single one on $k$-th position. Now, for each $k \in \mathbb N$ define $x_k \colon [0,1] \to l^2(\mathbb R)$ as the constant map taking the value $e_k$. Then, $\norm{x_k}_1 = \norm{e_k}_{l^2}=1$ and $\norm{C_f(x_k)}_1=\norm{f(e_k)}_{l^2}=\norm{ke_1}_{l^2}=k$ for $k \in \mathbb N$. This means that $C_f$ is not locally bounded.
\end{example}

We end the paper with a result providing necessary and sufficient conditions for a composition operator $C_f \colon BV_p([a,b],E) \to BV_q([a,b],E)$ to be locally bounded.

\begin{theorem}\label{thm:loc_bounded}
Let $E$ be a normed space and let $1\leq p \leq q<+\infty$. Assume that the composition operator $C_f$, generated by the mapping $f \colon E \to E$, maps the space $BV_p([a,b],E)$ into $BV_q([a,b],E)$. Then, the operator $C_f$ is locally bounded if and only if $f$ is H\"older continuous on bounded subsets of $E$ with exponent $p/q$.
\end{theorem}

\begin{proof}
First, we will show that $C_f$ is locally bounded. To this end let us fix $r>0$ and let the constant $L_r\geq 0$ be so that $\norm{f(u)-f(w)}_E\leq L_r \norm{u-w}^{p/q}_E$ for $u,w \in \ball_E(0,r)$. Set $R:=\norm{f(0)}_E+2L_r r^{p/q}$. Then, for every $x \in \ball_{BV_p}(0,r)$ and every $t \in [a,b]$ we have $\norm{x(t)}_E \leq \norm{x}_p\leq r$. So, $\var_q C_f(x) \leq L_r^q \var_p x$. Consequently, 
\begin{align*}
 \norm{C_f(x)}_q&= \norm[\big]{C_f(x)(a)}_E+(\var_q C_f(x))^{1/q} \leq \norm[\big]{f(0)}_E + \norm[\big]{f(x(a))-f(0)}_E + (\var_q C_f(x))^{1/q}\\
& \leq \norm[\big]{f(0)}_E +  L_r\norm[\big]{x(a)}^{p/q}_E + L_r(\var_p x)^{1/q} \leq \norm[\big]{f(0)}_E+ 2L_r r^{p/q} =R. 
\end{align*}
This shows that $C_f$ is locally bounded.

Now, let us assume that $C_f$ is locally bounded. Using a similar reasoning to the one we used in Example~\ref{ex:not_locally_bdd} we can show that then $f$ must be also locally bounded. Indeed, if there existed a sequence $(v_n)_{n \in \mathbb N}$ of elements of a bounded set $A$ with $\sup_{n \in \mathbb N}\norm{f(v_n)}_E=+\infty$, then for the sequence of constant maps $y_n(t)=v_n$ we would have $\sup_{n \in \mathbb N}\norm{C_f(y_n)}_q = \sup_{n \in \mathbb N}\norm{f(v_n)}_E = +\infty$. This would contradict the local boundedness of $C_f$. 

Suppose now that there is a closed ball $\ball_E(0,\rho)$ such that for every $n \in \mathbb N$ we can find two distinct points $u_n,w_n \in \ball_E(0,\rho)$ with $\norm{f(u_n)-f(w_n)}_E > 2Mn\norm{u_n-w_n}^{p/q}_E$, where $M:=\sup_{\norm{v}_E\leq \rho}\norm{f(v)}_E<+\infty$. Set $m_n:=[\norm{u_n-w_n}^{-p}_E]$; as in the proof of Theorem~\ref{thm:jospehy} by $[t]$ we denote the greatest integer smaller than or equal to $t$.  For each $n \in \mathbb N$ choose arbitrary points $a<a_1^n<a_2^n<\ldots<a_{m_n}^n<b$ and define
\[
 x_n(t)=\begin{cases}
			 u_n, & \text{if $t=a_k^n$ for some $k=1,\ldots,m_n$,}\\
			 w_n, & \text{if $t \in [a,b]\setminus\{a_1^n,\ldots,a_{m_n}^n\}$.}
			\end{cases}
\] 
It can be checked that $\var_p x_n \leq 2$, and so $\norm{x_n}_p \leq \rho+2^{1/p}$ for $n \in \mathbb N$.  But $\var_q C_f(x_n) \geq m_n\norm{f(u_n)-f(w_n)}_E^q \geq 2^q M^q n^q m_n \norm{u_n-w_n}^{p}_E\geq M^q n^q \to +\infty$ -- a contradiction. This ends the proof.
\end{proof}

\begin{remark}
The reason we do not consider bounded composition operators, that is, those whose range is a bounded subset of the codomain, is very simple. It turns out that any bounded composition operator $C_f$ mapping $BV_p([a,b],E)$ into $BV_q([a,b],E)$ for $1\leq p \leq q < +\infty$ is always constant. To prove this claim, suppose that $f(u)\neq f(0)$ for some non-zero $u \in E$. Further, fix $n \in \mathbb N$, take arbitrary $n$ points $a<a_1<a_2<\cdots<a_n<b$ and define $x \colon [a,b] \to E$ as
\[
 x_n(t)=\begin{cases}
         u, & \text{if $t=a_k$ for some $k=1,\ldots,n$,}\\
				 0, & \text{if $t \in [a,b] \setminus\{a_1,a_2,\ldots,a_n\}$.}
				\end{cases}
\]  
Then, $x_n \in BV_1([a,b],E)\subseteq BV_p([a,b],E)$ and $\norm{C_f(x_n)}_q \geq n^{1/q}\norm{f(u)-f(0)}_E$. This clearly means that $C_f$ is not bounded unless $f$ (and consequently, $C_f)$ is constant.
\end{remark}

\end{document}